\documentclass[12pt]{article}
\usepackage{amsmath,amstext, verbatim,remreset}
\usepackage[T2A]{fontenc}
\usepackage[cp1251]{inputenc}
\usepackage[english]{babel}
\usepackage{amsfonts,amssymb}
\usepackage{mathrsfs}
\newtheorem{Theorem}{Theorem}
\newtheorem{lemma}{Lemma}

\begin{document}
\begin{center}
{\bf \Large V. F. Babenko, O. V. Kovalenko}

O. Gonchar Dnipropetrovsk National University

babenko.vladislav@gmail.com

olegkovalenko90@gmail.com
\vskip 10mm

{\bf \Large On modulus of continuity of differentiation operator on weighted Sobolev classes}
\end{center}
\section {Introduction}

In the work of G.~Hardi and J.~Littlwood~\cite{HL} (1912) following theorem was proved (see also~\cite{BKKP}, Theorem~1.1.2).
\vskip 0.1 cm
{\bf Theorem A.} {\it Let a function $x(t)$ be defined for $t>0$ and its second derivative $x^{\prime\prime}(t)$ exist for $t>0$. Let, also, $f(t)$ and $g(t)$ be positive functions (both decreasing or both increasing). Then (as $t\to+\infty$) following statements hold:

$1.$ If $x(t)=O(f(t))$, $x^{\prime\prime}(t)=O(g(t))$, then $$x^{\prime}(t)=O\left(\sqrt{f(t)g(t)}\right);$$

$2.$ If $x(t)=o(f(t))$, $x^{\prime\prime}(t)=O(g(t))$, then $$x^{\prime}(t)=o\left(\sqrt{f(t)g(t)}\right);$$

$3.$ If $x(t)=O(f(t))$, $x^{\prime\prime}(t)=o(g(t))$, then $$x^{\prime}(t)=o\left(\sqrt{f(t)g(t)}\right).$$
}

This theorem had a great impact on formation of the whole field of problems connected with inequalities between derivatives. To confirm this, it is sufficient to note, that fundamental possibility of inequalities for upper bounds of derivatives of the functions, defined on the whole real line or half-line (see.~\cite{BKKP}, p.~18), can be easily attained from this theorem.

In 1928 L.~Mordell~\cite{Mordell} (see. also~\cite{BKKP}, Theorem~1.4.1) proved the following refinement of Theorem~A for non-increasing functions.  

{\bf Theorem~B.} 
{\it Let $f(t)$ and $g(t)$ be positive non-increasing on half-line ${\mathbb R}_+$ functions. 
If function $x(t)$ is defined on half-line ${\mathbb R}_+$ and for all $t>0$ there exists $x''(t)$ 
such that $$|x(t)|\leq f(t),\;|x''(t)|\leq g(t),$$ 
then for all $t>0$ $$|x^{\prime}(t)|\leq 2\sqrt{f(t)g(t)}.$$}

Let $I$ be a finite interval, whole real line ${\mathbb R}$, or positive half-line ${\mathbb R}_+$. Denote by $C^m(I)$ ($m \in {\mathbb Z}_+$) the set of all $m$ -- times continuously differentiable (continuous in the case $m=0$) functions $x:I \to {\mathbb R}$; by $L_\infty(I)$ we will denote the space of all measurable functions $x:I \to {\mathbb R}$ with finite norms
$$\|x\|_\infty:= {\rm ess\,sup}\,\{|x(t)|\,:\,t\in I\}.$$

Let $X$ be $C(I)$ or $L_\infty(I)$, $f\in C(I)$ is positive non-increasing function. For $x\in X$ set $$\|x\|_{X,f}:=\left\|\frac{x(\cdot)}{f(\cdot)}\right\|_X.$$

Using above notations, the result of the Theorem~B can be rewritten in the following way
\begin{equation}\label{vv1}
	\|x^\prime\|_{C({\mathbb R}_+),\sqrt{fg}}\leq 2\,\|x\|_{C({\mathbb R}_+),f}^{\frac 12}\| x^{\prime\prime}\|_{L_\infty({\mathbb R}_+),g}^{\frac 12}. 
\end{equation}

In the case, when $f(t)\equiv 1$ and $g(t)\equiv 1$ from~\eqref{vv1} we get Landau inequality~\cite{landau}, established in~1913:
\begin{equation}\label{vv2}
	\|x^\prime\|_{C({\mathbb R}_+)}\leq 2\,\|x\|_{C({\mathbb R}_+)}^{\frac 12}\| x^{\prime\prime}\|_{L_\infty({\mathbb R}_+)}^{\frac 12}\,. 
\end{equation}
Similar to~\eqref{vv2} sharp inequalities for functions, defined on the whole real line are also, in fact, contained in~\cite{landau} (see.~\cite{BKKP}, \S 1.2).

Later inequalities of the type~\eqref{vv2} for functions defined on ${\mathbb R}$ and ${\mathbb R}_+$ were generalized in many directions by many mathematicians. One of the brightest and most important results in the whole field is Kolmogorov inequality~\cite{Kolmogorov1,Kolmogorov2} for function defined on the real line ${\mathbb R}$. After this result the inequalities of the type~\eqref{vv2} are called Kolmogorov type inequalities. In articles~\cite{TihM,Argab,Arestov} and monographs~\cite{BKKP,Tihomirov,Fink} one can find detailed overview of known results about sharp inequalities for derivatives and further references.

We will discuss some of the results for functions defined on the half-line in a  more detailed way.

Let $T_r(t):=\cos r \arccos\,t,\,\,t\in[-1,1],$ be Chebyshev polynomials of first kind.
A.~P.~Matorin in 1955 proved the following theorem (see~\cite{Matorin}).
\vskip 0.1 cm
{\bf Theorem C. }{\it Let $k,r\in{\mathbb N},\,k<r.$ For arbitrary function $x\in L_{\infty,\infty}^r({\mathbb R}_+)$ the following inequality holds
	\begin{equation}\label{matorin}
		\|x^{(k)}\|_{\infty}\leq \frac {T_r^{(k)}(1)} {\left[T_r^{(r)}(1)\right]^{\frac kr}}\|x\|_\infty^{1-\frac kr}\|x^{(r)}\|_\infty^{\frac kr}.
	\end{equation}
	For the cases $r=2$ and $r=3$ inequality above is sharp.}

For $r>3$ inequality~\eqref{matorin} is not sharp. Sharp inequality, that estimates $\|x^{(k)}\|_{C({\mathbb R}_+)}$ using $\|x\|_{C({\mathbb R}_+)}$ and $\|x^{(r)}\|_{L_\infty({\mathbb R}_+)}$ for functions $x\in L_{\infty,\infty}^r({\mathbb R}_+)$ was received by I.~J.~Shoenberg and A.~Cavaretta (see.~\cite{Shenb1},~\cite{Shenb2}) in  1970 (see also~\cite{BKKP}, \S 3.3). 

For positive functions $f,g\in C(I)$ and natural $r$ set
$$L_{f,g}^r(I):=\left\{x\in C(I):\|x\|_{C(I),f}<\infty,\,x^{(r-1)}\in {\rm AC_{loc}},\,\|x^{(r)}\|_{L_\infty(I),g}<\infty\right\},$$
$$W_{f,g}^r(I):=\left\{x\in L_{f,g}^r(I)  :\|x^{(r)}\|_{L_\infty(I),g}\leq 1\right\}.$$

Function 
\begin{equation}\label{omega}
\omega(\delta)=\omega(D^k,\delta):=\sup\limits_{x\in W_{f,g}^r,\, \|x\|_{C(I),f}\leq\delta}\|x^{(k)}\|_{C(I)},\,\,\delta\geq 0
\end{equation}
is called modulus of continuity of $k$-th order differentiation operator on the class $W_{f,g}^r(I)$, ${(k=1,2,\dots,r-1)}$.

Note, that in the case $r = 2$ Theorem~B gives an estimate for the modulus of continuity $\omega(D^1,\delta)\leq 2\,\delta^\frac12$, $\delta > 0$. In the case of $f=g\equiv 1$ the result by I.~J.~Shoenberg and A.~Cavaretta  gives the value of $\omega(D^k,\delta)$ for all $\delta\geq 0$ (Theorem~C gives $\omega(D^k,\delta)$, $\delta\geq 0$ for $r\leq 3$).

Information about the connection between modulus of continuity of differentiation operator and Kolmogorov type inequalities and further references can be found in~\cite{BKKP}, \S1.7 and  Chapter~7.

The aim of this article is to study the function $\omega(D^k,\delta)$ for arbitrary $k,r\in{\mathbb N}$, $k<r$ and non-increasing continuous positive functions $f$ and $g$.

The article is organized in the following way. In \S 2 some auxiliary statements and in \S 3 --- main statements are given. \S 4 and \S 5 are devoted to proofs.

\section{Auxiliary results}

$\indent$ Let positive function $g \in C[a,b]$ be given. Function $G\in C^{r-1}[a,b]$ will be called perfect $g$-spline of the order $r$ with knots $a<t_1<\ldots<t_n<b$, if on each of the intervals $(t_i,t_{i+1}), i=0,1,\dots,n$ $(t_0:=a, t_{n+1}:=b)$ there exists derivative $G^{(r)}$ and $\frac{G^{(r)}(t)}{g(t)}\equiv\epsilon\cdot(-1)^i$ on intervals $(t_i,t_{i+1}), i=0,1,\dots,n$, where $\epsilon\in\{1,-1\}$.

Denote by $\Gamma^{r}_{n,g}[0,a]$ the set of all perfect $g$-splines $G$, defined on $[0,a]$, of the order $r$ with not more than $n$ knots.

Below $f,g$ will denote continuous positive non-increasing on $[0,\infty)$ functions.

The next theorem proofs existence and some properties of perfect $g$-spline $G_{r,n,f,a}\in \Gamma_{n,g}^r[0,a]$ least deviating from zero in $\|\cdot\|_{C[0,a],f}$  norm.

\begin{Theorem}\label{th::1}   Let numbers $a>0$, $r\in {\mathbb N}$, $n\in {\mathbb Z}_+$ be given. Then there exists a perfect $g$-spline $G_{r,n,f,a}\in \Gamma^{r}_{n,g}[0,a]$, that has $n+r+1$ oscillation points, i.~e. such, that there exist $n+r+1$ points $0\leq t_1<t_2<\ldots<t_{n+r+1}= a$, such that
\begin{equation}\label{th1.0}
	G_{r,n,f,a}(t_i)=(-1)^{i+r+1}\|G_{r,n,f,a}\|_{C[0,a],f}\cdot f(t_i),\, i=1,2,\dots,n+r+1.
\end{equation}
Perfect $g$-spline, that satisfies condition~\eqref{th1.0} is unique in the set $\Gamma^{r}_{n,g}[0,a]$.

For $a>0$ set $\varphi_{r,n,f}(a):= \|G_{r,n,f,a}\|_{C[0,a],f}$. Then function $\varphi_{r,n,f}(a)$ continuously depends on $a\in(0,\infty)$ and increases with growing of $a$.
\end{Theorem} 

The role of perfect $g$-spline becomes clearer due to the following theorem.

\begin{Theorem}\label{th::3} 
	Let  $r\in {\mathbb N}$, $n\in {\mathbb Z}_+$ and $\delta>0$ be such, that {${\varphi_{r,n,f}(a)=\delta}$} for some $a>0$. Then for $k=1,2,\dots,r-1$ $$\omega(D^k,\delta)\leq \left\|G_{r,n,f,a}^{(k)}\right\|_{C[0,a]}.$$
\end{Theorem}
\section {Main results}

If $f(t)\equiv 1$ and  $g(t)\equiv 1$, then $\varphi_{r,n,f}(\infty):=\lim\limits_{a\to+\infty} \varphi_{r,n,f}(a)=\infty$ for all $r\in{\mathbb N}$, $n\in{\mathbb Z}_+$. In the case, when $f,g$ are arbitrary positive non-increasing continuous functions this is not always true.
 
Set $g_k(t):=\int\limits_0^tg_{k-1}(s)ds,\,k=1,2,\dots,{r},$ where $g_0:=g$. The following theorem holds.

\begin{Theorem}\label{th::2}
	Let numbers $n\in {\mathbb Z}_+$ and $r\in{\mathbb N}$ be given. $\varphi_{r,n,f}(\infty)<\infty$ if and only if, the following conditions hold
	$$A_0:=\int\limits_0^\infty g(t)dt<\infty,$$
	\begin{equation} \label {th1}
		A_k:=\int\limits_0^\infty\left[\sum\limits_{s=0}^{k-1}\frac{(-1)^{k-s-1}A_s}{(k-s-1)!}t^{k-s-1}+(-1)^kg_k(t)\right]dt<\infty,\,k=1,\dots,r-1,
	\end{equation}
	\noindent and
	\begin{equation} \label {th2}
		\sup\limits_{t\in[0,\infty)}\frac{\left|\sum\limits_{s=0}^{r-1}\frac{(-1)^{r-s-1}A_s}{(r-s-1)!}\,t^{r-s-1} +(-1)^{r} g_{r}(t)\right|}{f(t)}<\infty. 
	\end{equation}
\end{Theorem}

{\bf Remark. }{\it From Theorem~\ref{th::2} if follows that for all $r\in{\mathbb N}$, $n\in{\mathbb Z}_+$ $\varphi_{r,n,f}(\infty)<\infty$ if and only if $\varphi_{r,0,f}(\infty)<\infty$.}

In the case, when conditions~\eqref{th1} hold, set
\begin{equation}\label{Pk}
	P_k(t) := \sum\limits_{s=0}^{k-1}\frac{(-1)^{k-s-1}A_s}{(k-s-1)!}t^{k-s-1}+(-1)^kg_k(t),\, k=1,2,\dots, r.
\end{equation}

If $\varphi_{r,0,f}(\infty)=\infty$,  then, in virtue of Theorems~\ref{th::1} and~\ref{th::2}, for all $r\in{\mathbb N}$, $n\in{\mathbb Z}_+$ and  $ \delta>0$ there exists a number ${\delta_{r,n}>0}$ such that $\|G_{r,n,f,\delta_{r,n}}\|_{C[0,\delta_{r,n}],f}=\delta.$ In this case for all $\delta>0$ the function $\omega(D^k,\delta)$ is characterized by the following theorem.

\begin{Theorem}\label{th::4} 
	Let $r\in{\mathbb N}$ and $\varphi_{r,0,f}(\infty)=\infty$. Then for all $\delta>0$ and $k=1,2,\dots,r-1$ $$\omega(D^k,\delta)=\lim\limits_{n\to\infty}\left|G^{(k)}_{r,n,f,\delta_{r,n}}(0)\right|.$$
\end{Theorem}
Information about the function $\omega(D^k,\delta)$ in the case, when $\varphi_{r,0,f}(\infty)<\infty$ is given by the following theorem.
\begin{Theorem}\label{th::5} 
	Let $r\in{\mathbb N}$, $n\in{\mathbb Z}_+$ and $\varphi_{r,0,f}(\infty)<\infty$.
	Then for all $k=1,2,\dots,r-1$ 
	$$\omega(D^k,\varphi_{r,n,f}(\infty))=\lim\limits_{a\to\infty}\left|G^{(k)}_{r,n,f,a}(0) \right|.$$
\end{Theorem}

In the case, when $\varphi_{r,0,f}(\infty)<\infty$, information about asymptotic behavior of the function $\varphi_{r,n,f}(\infty)$ as $n\to\infty$ and fixed $r$  is given by the following theorem.

\begin{Theorem}\label{th::6}
	Let $r\in{\mathbb N}$ and $\varphi_{r,0,f}(\infty)<\infty$. $\lim\limits_{n\to\infty} \varphi_{r,n,f}(\infty)>0$ if and only if $\varliminf\limits_{t\to\infty}\frac{f(t)}{\left|P_r(t)\right|\vphantom{^2}}<\infty$, 
	where the function $P_r(t)$ is defined in~\eqref{Pk}. 
\end{Theorem}

\section{Proof of auxiliary results}
\subsection{Proof of Theorem~\ref{th::1}.}
$\indent$ Proof of existence and uniqueness of perfect $g$-spline $G_{r,n,f,a}$ uses ideas that were used to prove Theorem~3.3.1 in monography~\cite{BKKP}. 
 
In the space ${\mathbb R}^{n+1}_1$ consider sphere $S^n$ with radius $a$, i.~e. $$S^n=\left\{ \xi=(\xi_1,\xi_2,\dots,\xi_{n+1}):\sum_{i=1}^{n+1}|\xi_i|=a\right\}.$$ For each $\xi\in S^{n}$ consider partition of the segment $[0,a]$ by points $$\xi_0:=0,|\xi_1|,|\xi_1|+|\xi_2|,\dots,\sum_{i=1}^{n}|\xi_i|,\sum_{i=1}^{n+1}|\xi_i|=a.$$
Set $I_k:=\left(\sum_{i=0}^{k-1}|\xi_i|,\sum_{i=0}^{k}|\xi_i|\right)$, $k=1,2,\dots,n+1$. For each of partitions consider a function $$G_{\xi}^a(t):=\frac 1{(r-1)!}\int_0^a{(t-u)^{r-1}_+g^a_{\xi}(u)du},$$ where $g^a_{\xi}(t)=g(t)\,{\rm sgn}\xi_{k}$ on each segment $I_k$, $(k=1,2,\dots,n+1).$

Then we have $(G_{\xi}^a)^{(r)}=g_{\xi}^a$, and hence $G_{\xi}^a$ is $g$-spline with knots in the points of partition. Let $Q^{\xi,a}_{n+r-1}(t)=\sum_{i=0}^{n+r-1}{a_i(\xi)t^i}$ be the polynomial on which $\inf\limits^{}_{Q_{n+r-1}}\|G_{\xi}^{a}-Q_{n+r-1}\|_{C,f[0,a]}$ is attained. Consider mapping $\phi:S^n\to{\mathbb R}^n$, $\phi(\xi):=(a_r(\xi),\dots,a_{n+r-1}(\xi))$. From the definition of $\phi$ and properties of polynomials of the best approximation it follows that $\phi$ is continuous and odd. Hence from Borsuk's theorem it follows that there exists $\xi_0\in S^n$, such that $\phi(\xi_0)=0$. This means that the polynomial $Q^{\xi_0,a}_{n+r-1}$ has order not bigger than $r-1$. Hence for the function $G_{r,n,f,a}:=G_{\xi_0}-Q^{\xi_0,a}_{n+r-1}$ we have $G_{r,n,f,a}^{(r)}=g_{\xi_0}$. Hence, $G_{r,n,f,a}$ is $g$-spline. In virtue of generalization of Chebyshev's theorem about oscillation (see, for example,~\cite{KK}, Chapter~9, \S5) $G_{r,n,f,a}$ has $n+r+1$ oscillation points $0\leq t_1<t_2<\ldots<t_{n+r+1}\leq a$, and hence at least $n+r$ sign changes. Hence in virtue of Rolle's theorem $G_{r,n,f,a}^{(r)}$ has $n$ sign changes. This means that $G_{r,n,f,a}\in \Gamma_{n,g}^{r}$. 

We will prove that $t_{n+r+1}=a$. In the opposite case taking into account that $f$ is non-increasing, we would receive that $G_{r,n,f,a}^{'}$ has $n+r$ zeroes and hence $G_{r,n,f,a}^{(r)}$ would have $n+1$ sign changes, which is impossible. Multiplying, if needed, the function $G_{r,n,f,a}$ by $-1$, we will get perfect $g-$spline, for which the equalities~\eqref{th1.0} hold.

To prove uniqueness it is sufficient to prove the following lemma.
\begin{lemma}\label{l::pr3} 
	Let $a>0$, $r\in{\mathbb N}$, $n\in{\mathbb Z}_+$ be given. Let $G_{r,n,f,a}$ be the perfect $g$-spline,  for which equalities~\eqref{th1.0} hold. Let $\left\{u_k\right\}_{k=1}^n$ be its knots. Let also $s\in{\mathbb N}$ and spline $G\neq \pm G_{r,n,f,a}\in \Gamma_{s,g}^r[0,a]$ such, that $\|G\|_{C[0,a],f}\leq \|G_{r,n,f,a}\|_{C[0,a],f}$  be given. Let $\left\{v_k\right\}_{k=1}^s$ be the knots of spline $G$. Then  $s>n$ and, if $s = n+1$, then $u_i>v_i$, $i=1,\dots,n$.
\end{lemma}

We can count, that ${\rm sgn}G^{(r)}(0) = {\rm sgn}G_{r,n,f,a}^{(r)}(0)$. Set $\delta(t):=G_{r,n,f,a}(t)-G(t) \,,\,t\in[0,a]$. $\delta(t)$ has $n+r$ zeroes on $[0,a]$, hence $\delta^{(r)}(t)$ has $n$ sign changes in $(0,a)$. From the other side, $\delta^{(r)}(t)$ can't have sign changes inside of the intervals $(u_i,u_{i+1})$ $(i=0,\dots,n)$, $(u_0:=0, u_{n+1}:=a)$. Thus, function $\delta^{(r)}(t)$ has exactly $n$ sign changes, at that it is not equal to $0$ identically in any of the intervals $(u_i,u_{i+1})$ $(i=0,\dots,n)$. This means, that on each interval $(u_i,u_{i+1})$ $(i=0,\dots,n)$ there exists at least one knot of the spline $G$, and hence $s > n$. If $s = n+1$, then on each interval $(u_i,u_{i+1})$ $(i=0,\dots,n)$ there exists exactly one knot of the spline $G$. Lemma is proved.

For fixed $r\in {\mathbb N}$ and $n\in{\mathbb Z}_+$ $\varphi_{r,n,f}(a)$ increases when $a$ grows in virtue of Lemma~\ref{l::pr3}. Its  continuity follows from continuity of functions $f$ and $g$. Theorem is proved.

\subsection {Proof of Theorem~\ref{th::3}.}

We will need the following lemma.
\begin{lemma}\label{l::signChange}
Let $r\in{\mathbb N}$, $n\in{\mathbb Z}_+$, $a>0$ and $x\in  L^r_{\infty,\infty}[0,a]$ be given. Assume function $x$ has at least $n+r$ sign changes and function $x^{(r)}$ has not more than $n$ sign changes. Then for all $s=0,1,\dots,r-1$ 
\begin{equation}\label{th3.0}
{\rm sgn }x^{(s)}(0) = -{\rm sgn }x^{(s+1)}(0).
\end{equation}
\end{lemma}
{\bf Remark.} Notation ${\rm sgn }x^{(s)}(0) = \pm 1$ means that there exists $\varepsilon > 0$ such, that ${\rm sgn }x^{(s)}(t)=\pm 1$  in neighborhood $(0,\varepsilon)$ (almost everywhere in neighborhood $(0,\varepsilon)$ in the case $s=r$). 

From conditions of the lemma it follows, that the function $x^{(s)}$ has exactly $n+r-s$ sign changes, $s = 0,1,\dots,r$. If we assume that for some $0\leq s \leq r-1$ equality~\eqref{th3.0} does not hold, then in virtue of Fermat's theorem we will get extra sign change of the function $x^{(s+1)}$. Lemma is proved.

Let's return to the proof of the theorem. Assume contrary, let a function $x\in W^r_{f,g}({\mathbb R}_+)$ be such, that $\|x\|_{C[0,\infty),f}\leq \delta$ and  $\|x^{(k)}\|_{C[0,\infty)}> \left\|G_{r,n,f,a}^{(k)}\right\|_{C[0,a)}$. We can count, that $\left|x^{(k)}(0)\right|>\left|G_{r,n,f,a}^{(k)}(0)\right|$ (if this is not true, then there exists a point $t_0>0$ such, that $|x^{(k)}(t_0)|>\left\|G_{r,n,f,a}^{(k)}\right\|_{C[0,a)}$  and we can instead of $x(t)$ consider the function $y(t):= x(t+t_0)$. At that we will get  $y\in W^r_{f,g}({\mathbb R}_+)$ and $\|y\|_{C[0,\infty),f}\leq \delta$, and uniform norms of the function $x$ and its derivatives will not increase, because  $f$ and $g$ are non-increasing). Moreover, multiplying, if needed, functions $x$ and $G_{r,n,f,a}$ by $-1$ we can count that
\begin{equation}\label{th3.00}
x^{(k)}(0)>G_{r,n,f,a}^{(k)}(0)>0.
\end{equation}

Set $\Delta(t):=x(t)-G_{r,n,f,a}(t)$. Note, that in virtue of building of $g$-splines $G_{r,n,f,a}(t)$, functions $\Delta(t)$ and $G_{r,n,f,a}(t)$ have not less than $n+r$ sign changes (function $G_{r,n,f,a}(t)$ has exactly $n+r$ sign changes), and functions $\Delta^{(r)}(t)$ and $G^{(r)}_{r,n,f,a}(t)$ can have sign changes only in the knots of $g$-spline $G_{r,n,f,a}(t)$, and hence not more than $n$ sign changes. From Lemma~\ref{l::signChange} and~\eqref{th3.00} we get
\begin{equation}\label{th3.02}
(-1)^k G_{r,n,f,a}(0)>0.
\end{equation}

In virtue of~\eqref{th3.02} $(-1)^k\Delta(0) < 0$, and hence in virtue of Lemma~\ref{l::signChange} we get $\Delta^{(k)}(0) < 0$. But this contradicts to~\eqref{th3.00}. Theorem is proved.  

\section {Proof of main results}

\subsection {Proof of Theorem~\ref{th::2}.}

We will prove first, that the statement of the theorem is true in the case $n=0$. In the case $n=0$ we will write $\varphi_{r,f}$ instead of $\varphi_{r,0,f}$ and $G_{r,f,M}$ instead of $G_{r,0,f,M}$. 
To prove the theorem we will need following lemma.
\begin{lemma}\label{l::01}
	Let condition~\eqref{th1} hold. Let $M>0$ and $h_m(t)$ is $m$-th primitive of the function $g(t)$ on interval $[0,M]$, which has $m$ zeroes $(1\leq m\leq r)$. Denote by $\alpha_m$ the first zero of the function $h_m(t)$. 	Then the following inequalities hold
	$$|h_m(t)|<|P_m(t)|\,,t\in[0,\alpha_m]$$
	and $$|P_m(t)-P_m(0)+h_m(0)|\leq|h_m(t)|\,,t\in[0,\gamma_m],$$
	where $\gamma_m$ is zero of the function $P_m(t)-P_m(0)+h_m(0)$.
\end{lemma} 

We will prove the first inequality, the second can be proved using similar arguments.
We will proceed using induction on $m$. 

Let $m=1$. $P_1(t)=-A_0+g_1(t)$. Since condition~\eqref{th1} holds and $P_1(0)=-A_0$, we have $P_1(\infty)=0$. Since the function $h_1(t)=g_1(t)+C$ has one zero, $C\in (-A_0,0]$, this means $|h_1(t)|<|P_1(t)|\,\forall t\in[0,\alpha_1]$.

Let the statement of the lemma hold in the case $m=k\leq r-1$. We will show that it holds in the case $m=k+1$ too. Let the number $k+1$ be even. Then $P_{k+1}(0)>0$ and, in virtue of Lemma~\ref{l::signChange},
\begin{equation}\label{l1.1}
	h_{k+1}(0)>0.
\end{equation}
 
Assume contrary. Let a point $t_0\in [0,\alpha_{k+1})$ such that $h_{k+1}(t_0)\geq P_{k+1}(t_0)$ exist. Since the function $h_{k+1}(t)$ has $k+1$ zeroes, the function $h^{(k+1)}_{k+1}(t)= g(t)$ does not have zeroes, we have that the function $h_{k+1}^{'}(t)$ has $k$ zeroes, and hence in virtue of induction assumption and $P_{k+1}^{'}(t)=P_k(t)$ we get that $|h_{k+1}^{'}(t)|<|P_k(t)|,\,\forall t\in [0,\beta]$, where $\beta$ --- is first zero of the function $h_{k+1}^{'}(t)$. That's why in virtue of~\eqref{l1.1} for all $t\in[0,\beta]$ the following inequality holds
\begin{equation}\label{l1.2}
	0\geq h_{k+1}^{'}(t)>P_k(t). 
\end{equation}
In virtue of Rolle's theorem $\beta>\alpha_k$, which means, that inequality~\eqref{l1.2} is valid for all $t\in[0,\alpha_k]$.
Since
$h_{k+1}(\alpha_{k+1})=0<P_{k+1}(\alpha_{k+1})$, we have $P_{k+1}(t_0)-P_{k+1}(\alpha_{k+1})<h_{k+1}(t_0)-h_{k+1}(\alpha_{k+1})$, from the other side inequality~\eqref{l1.2} holds, and hence
$$P_{k+1}(t_0)-P_{k+1}(\alpha_{k+1})=-\int\limits_{t_0}^{\alpha_{k+1}}P_k(t)dt>-\int\limits_{t_0}^{\alpha_{k+1}}{h_{k+1}^{'}
(t)dt}= $$ $$=h_{k+1}(t_0)-h_{k+1}(\alpha_{k+1}).$$
Contradiction. Lemma is proved.

\vskip 0.2 cm
Let's return to the proof of the theorem in the case when $n=0$. 

Let conditions~\eqref{th1} and~\eqref{th2} hold. Set $K_r:=\sup\limits_{t\in[0,\infty)}\frac{|P_r(t)|}{f(t)}$. Assume contrary, let ${\varphi_{r,f}(\infty)=\infty}$. This means that there exists $M>0$ such that $\varphi_{r,f}(M)>K_{r}$. In virtue of the Lemma~\ref{l::01} on interval $[0,\alpha_{r}^M]$ the following inequality $|G_{r,f,M}(t)|<|P_{r}(t)|$  holds. Moreover, $|P_{r}(t)|\leq K_r\,f(t)<\varphi_{r,f}(M)\,f(t)$. Thus, on interval $[0,\alpha_{r}^M]$ the following inequality $|G_{r,f,M}(t)|<\varphi_{r,f}(M)\,f(t)$ holds. However in this case $G_{r,f,M}(t)$ has not more than $r$ oscillating points --- contradiction. Sufficiency is proved.

Let now $\varphi_{r,f}(\infty)< \infty$. This means that for all $a, t\geq 0$ we have $\left|G_{r,f,a}(t)\right| \leq \varphi_{r,f}(\infty) f(t) \leq \varphi_{r,f}(\infty) f(0)$. Passing to the limit, when $a\to\infty$ we get existence of bounded on $[0,\infty)$ primitive $Q_r$ of order $r$ of the function $g(t)$. Since functions $f(t)$ and $g(t)$ are bounded, we get that all functions $Q_r^{(k)}(t)$, $k=1,\dots,r-1$ are also bounded on $[0,\infty)$. Note, that the only bounded on $[0,\infty)$ primitives of the function $g(t)$ of order $k\in{\mathbb N}$ are functions $P_k(t) + C_k$, where $C_k\in{\mathbb R}$, and only in the case when corresponding conditions~\eqref{th1} hold.  This means that conditions~\eqref{th1} hold. Necessity of conditions~\eqref{th1} is proved.

Note, that from arguments above it follows that the following lemma holds

\begin{lemma}\label{l::02}
	Let $r\in{\mathbb N}$, $r\geq 2$ and $\varphi_{r,f}(\infty)<\infty$. Then $\left|G_{r,f,M}^{(r-k-1)}(0)\right|\to A_k$ and $\alpha_{k+1}^M\to\infty$ when $M\to \infty$, where $\alpha_{k+1}^M$ is the first zero of the function 
	$G_{r,f,M}^{(r-k-1)}$ $(k=0,1,\dots,r-2)$.
\end{lemma} 

We will prove, that condition~\eqref{th2} also holds. If $f(\infty)>0$, then condition~\eqref{th2} holds always when conditions~\eqref{th1} hold. Below we will count that
\begin{equation}\label{th6}
f(\infty)=0.
\end{equation}
From Lemma~\ref{l::01} we get

$$\left|\int\limits_0^{\alpha_{r-1}^M}P_{r-1,f}(t)dt\right|\geq\left|G_{r,f,{M}}(0)-G_{r,f,M}(\alpha_{r-1}^{M})\right| \geq$$
$$  \left|\int\limits_0^{\gamma_M}P_{r-1,f}(t)-P_{r-1}(0)+G_{r,M}'(0)dt\right|.$$

In virtue of the Lemma~\ref{l::02} and equality~\eqref{th6} we  get 
\begin{equation}\label{th7}
	\left|\,G_{r,f,{M}}(0)\right|\to A_{r-1}. 
\end{equation}
when $M\to\infty$.

We will show that $\varphi_{r,f}(\infty)\geq \sup\limits_{t\in[0,\infty)}\frac{P_r(t)}{f(t)}$. Assume contrary. Let a point $t_0$ such, that $\left|P_{r}(t_0)\right|>\varphi_{r,f}(\infty)f(t_0)$, exists. We can choose $\varepsilon>0$ in such a way, that
\begin{equation}\label{th3}
\left|P_{r}(t_0)-{\rm \varepsilon\,sgn\left[P_{r}(0)\right]} \right|>\varphi_{r,f}(\infty)f(t_0).
\end{equation}

In virtue of~\eqref{th7} and Lemma~\ref{l::01} we can choose $M>0$ enough big, so that $$\left|P_{r}(t)-{\rm \varepsilon}\,sgn\left[P_{r}(0)\right] \right|<\left|G_{r,f,M}(t)\right|<\left|P_{r}(t)\right|\,\forall t\in[0,\gamma],$$ where $\gamma$ is zero of the function $P_{r}(t)-{\rm \varepsilon}\,sgn\left[P_{r}(0)\right].$ Since inequality~\eqref{th3} holds, we have $\gamma>t_0$, and hence $\left|G_{r,f,M}(t_0)\right|>\varphi_{r,f}(\infty)f(t_0)$ --- contradiction.
Thus condition~\eqref{th2} is proved. Theorem is proved in the case, when $n=0$.

Let $n$ be arbitrary natural number now.

We will prove, that for all $r,n\in{\mathbb N}$ $\varphi_{r,n,f}(\infty)<\infty$ if and only if $\varphi_{r,f}(\infty)<\infty$.

It is clear, that for all $M>0$ $\varphi_{r,f}(M)\geq\varphi_{r,n,f}(M)$, and hence $\varphi_{r,f}(\infty)<\infty$ implies $\varphi_{r,n,f}(\infty)<\infty$.

Assume $\varphi_{r,n,f}(\infty)<\infty$. Denote by $t_{n,k}^M$ $k$-th knot of $g$-spline $G_{r,n,f,M}(t)$, $k=1,2,\dots,n$. 
Set $t_{n,0}^M:=0$, $t_{n,n+1}^{M}:=M$. Let $1\leq k\leq n+1$ be the smallest number of the knots of $g$-spline $G_{r,n,f,M}(t)$, for which the set $\left\{t_{n,k}^M\colon M>0\right\}$ is unbounded. We can choose increasing sequence $\left\{M_l\right\}_{l=1}^{\infty}$, $M_l\to\infty$ as $l\to\infty$ such, that $t_{n,s}^{M_l}\to t_{n,s}<\infty$,  $s\leq k-1$ and $t_{n,k}^{M_l}\to \infty$, as $l\to\infty$. 

Denote by $G_{r,f,M}^K(t)$ the least deviation from zero in norm $\|\cdot\|_{C[K,K+M],f}$ primitive of the order $r$ of the function $g$ on the segment $[K,K+M]$, set $$\varphi_{r,f}^K(M):=\|G_{r,f,M}^K\|_{C[K,K+M],f}.$$ Then for all $l$ $\|G_{r,n,f,M_l}\|_{C[t_{n,k-1}^{M_l},t_{n,k}^{M_l}],f}\geq \varphi_{r,f}^{t_{n,k-1}}(t_{n,k}^{M_l}-t_{n,k-1}^{M_l})$. Passing to the limit as $l\to\infty$ we get $\varphi^{t_{n,k-1}}_{r,f}(\infty)\leq \varphi_{r,n,f}(\infty)<\infty$. To finish the proof of the theorem is is sufficient to note, that from proved above case when $n=0$ it follows, that for all $K>0$ $\varphi_{r,f}^K(\infty)<\infty$ if and only if $\varphi_{r,f}(\infty)<\infty$. Theorem is proved.

\vskip 0.1 cm
{\bf Remark. } In the case, when $f\equiv 1$, for all natural $r$ $\varphi_{r,f}(\infty)<\infty$ implies $\varphi_{r-1,f}(\infty)<\infty$. At the same time not for all functions $f$ $\varphi_{r,f}(\infty)<\infty$ implies $\varphi_{r-1,f}(\infty)<\infty$.
\vskip 0.1 cm
Really, in the case, when $f\equiv 1$ condition~\eqref{th2} holds always, when condition~\eqref{th1} holds. If $$g(t)=\left[-1+(t+\sqrt{3})^2\right]e^{-\frac{(t+\sqrt{3})^2}{2}}$$ and $$f(t)=e^{-\frac{(t+\sqrt{3})^2}{2}},$$ then $$P_1(t)=-\left[(t+\sqrt{3})\right]e^{-\frac{(t+\sqrt{3})^2}{2}},$$ $$P_2(t)=e^{-\frac{(t+\sqrt{3})^2}{2}}$$ and condition~\eqref{th2} holds when $r=2$ and does not hold when $r=1$.

\subsection {Proof of the Theorem~\ref{th::4}.}
To prove Theorem~\ref{th::4} it is sufficient to prove that for all $\delta>0$ there exists defined on $[0,\infty)$ $g$-spline $G_{r,\delta}=G_{r,\delta}(\cdot,f,g)$ of order $r$ with infinite number of knots $y_k$, $(k=1,2,\dots)$ $0=y_0<y_1<\ldots<y_k<\ldots$, $y_k\to\infty\,(k\to\infty)$, with following properties:

$1.$ $\|G_{r,\delta}\|_{C[0,\infty),f}= \delta$ and $G_{r,\delta}^{(r)}\equiv g$ or $G_{r,\delta}^{(r)}\equiv -g$ on intervals $(y_k,y_{k+1})$ $(k=0,1,2,\dots)$.

$2.$ For all $c>0$ sequences $\left\{ G_{r,n,f,\delta_{r,n}}^{(k)}\right\}_{n=0}^{\infty}$ $(k=0,1,\dots,r-1)$ (whose elements are defined on $[0,c]$ for big enough $n$) converge to $G_{r,\delta}^{(k)}$ uniformly on $[0,c]$.

Really, from Theorem~\ref{th::3} it will follow, that $\omega(D^k,\delta)=\left|G^{(k)}_{r,\delta}(0)\right|,$ and from condition~2 it will follow, that $$\lim\limits_{n\to\infty}\left|G^{(k)}_{r,n,f,\delta_{r,n}}(0)\right|=\left|G^{(k)}_{r,\delta}(0)\right|.$$ 

Lemma~\ref{l::pr3} implies, that the sequence $\left\{\delta_{r,n}\right\}_{n=0}^\infty$ strictly increases.
Moreover, this sequence is unbounded, because otherwise we would get perfect $g$-spline $G$ with arbitrarily close oscillating points, which is impossible because functions $G$ and $G^{(r)}$ (and hence $G^{\,'}$) are bounded. 

Denote by $t_{n,k}$ $(k=1,\dots,n,\,n=1,2,\dots)$ the knots of the $g$-spline $G_{r,n,f,\delta_{r,n}}$. Lemma~\ref{l::pr3} implies, that every sequence $\left\{t_{n,k}\right\}_{n=k}^\infty\, (k=1,2,\dots)$ decreases, and hence has a limit.

Let $0\leq y_1<y_2<\ldots $ be all distinct finite limits of these sequences, ordered in ascending way. Then for all $i\in{\mathbb N}$ and for all small enough $\varepsilon>0$ there exists $N=N(i,\varepsilon)$ such, that for all $n>N(i,\varepsilon)$ $G_{r,n,f,\delta_{r,n}}^{(r)}\equiv g$ or $G_{r,n,f,\delta_{r,n}}^{(r)}\equiv -g$ on $I_i(\varepsilon):=(y_{i-1}+\varepsilon,y_{i}-\varepsilon)$. Since $\varepsilon>0$ is arbitrary, we get existence of point-wise limit $\lim\limits_{n\to\infty}G_{r,n,f,\delta_{r,n}}=:G_{r,\delta}$, at that on intervals $(y_i,y_{i+1})$ $G_{r,\delta}^{(r)}\equiv g$ or $G_{r,\delta}^{(r)}\equiv -g$ $(i=1,2,\dots)$. It is clear, that $\|G_{r,\delta}\|_{C[0,\infty),f}= \delta$. Since $\lim\limits_{n\to\infty}{\delta_{r,n}}=+\infty$, then $y_k\to\infty\,(k\to\infty)$. 

Let's fix some $c>0$. Starting with some $n$ all $g$-splines $G_{r,n,f,\delta_{r,n}}$ are defined on $[0,c]$. From $\|G_{r,n,f,\delta_{r,n}}\|_{C[0,a_{r,n}],f}=\delta$ and the fact that $f$ is non-increasing is follows that the sequence  $\left\{G_{r,n,f,\delta_{r,n}}\right\}_{n=0}^{\infty}$ is uniformly bounded; from $\left|G_{r,n,f,\delta_{r,n}}^{(r)}(t)\right|\leq g(t)$ almost everywhere on $[0,\infty)$ and the fact that $g$ is non-increasing it follows that sequences  $\left\{G_{r,n,f,\delta_{r,n}}^{(k)}\right\}_{n=0}^{\infty}$, $k=0,\dots,r-1$, are uniformly bounded and equicontinuous. The later implies uniform convergence on $[0,c]$ of the sequence $G_{r,n,f,\delta_{r,n}}$ to $G_{r,f}$.  Theorem is proved.

\subsection{Proof of the Theorem~\ref{th::5}.} 

Let $n\geq 0$. We can choose increasing sequence $\left\{M_k\right\}_{k=1}^{\infty}$, $M_k\to\infty$ as  $k\to\infty$ in such a way, that all sequences $t_{n,s}^{M_k}$, $1\leq s\leq n$ (as above, $t_{n,s}^{M_k}$ is $s$-th knot of $g$-spline $G_{r,n,f,M_k}$) have limits (finite or infinite).
Let $t_{n,1}<\ldots<t_{n,m}$ be all distinct finite limits of these sequences in ascending order. Analogously to the proof of Theorem~\ref{th::4} we get uniform on each segment $[0,c],\,c>0$ convergence of the sequence $G_{r,n,f,M_k}$ to $g$-spline $P_{r,n,f,\{M_k\}}$ with $m$ knots (defined on the whole half-line) together with all derivatives up to the order $r-1$ inclusively. For brevity we will write $P_{r,n,f}$ instead of $P_{r,n,f,\{M_k\}}$.

Let function $x(t)$ be such, that
\begin{eqnarray}\label{th3.1}
\|x\|_{C[0,\infty),f}\leq\varphi_{r,n,f}(\infty),
\nonumber\\ \left\|x^{(r)}\right\|_{L_\infty(0,\infty),g}\leq 1.
\end{eqnarray}
We will show that for all $s=1,2,\dots,r-1$ 
\begin{equation}\label{ekstr}
	\left\|x^{(s)}\right\|_{C[0,\infty)}\leq \left\|P^{(s)}_{r,n,f}\right\|_{C[0,\infty)}.
\end{equation}

Assume contrary, let for some $s$ $\left\|x^{(s)}\right\|_{C[0,\infty)}> \left\|P^{(s)}_{r,n,f}\right\|_{C[0,\infty)}.$ Then there exists $\varepsilon>0$ such, that $\left\|x^{(s)}\right\|_{C[0,\infty)}> (1+\varepsilon)\left\|P^{(s)}_{r,n,f}\right\|_{C[0,\infty)}.$ We can count, that $\left|x^{(s)}(0)\right|>(1+\varepsilon)\left|P_{r,n,f}^{(s)}(0)\right|$ (if this is not true, then there exists a point $t_0>0$ such, that $\left|x^{(s)}(t_0)\right|>\left\|(1+\varepsilon)P_{r,n,f}^{(s)}\right\|_{C[0,\infty)}$ and instead of the function $x(t)$ we can consider $y(t):= x(t+t_0)$. At that since the functions $f$ and $g$  are non-increasing, the conditions~\eqref{th3.1} are not broken and uniform norms of the function $x$  and its derivatives do not increase).
Moreover, we can count that
\begin{equation}
	x^{(s)}(0)>(1+\varepsilon)P_{r,n,f}^{(s)}(0)>0
\end{equation}
(if this is not true, we can multiply $x$ and (or) $P_{r,n,f}$ by $-1$).
Set $\Delta_k(t):=x(t)-(1+\varepsilon)P_{r,n,f,M_k}(t)$, $t\in [0,M_k]$. We can choose $k$ so big, that
\begin{equation}\label{th4.2}
	x^{(s)}(0)>(1+\varepsilon)P^{(s)}_{r,n,f,M_k}(0)
\end{equation}
and
\begin{equation}\label{th3.3}
	(1+\varepsilon)\varphi_{r,n,f}(M_k)>\varphi_{r,n,f}(\infty).
\end{equation}

From Lemma~\ref{l::signChange} we get
\begin{equation}\label{th4.02}
	(-1)^s P_{r,n,f,M_k}(0)>0.
\end{equation}

In virtue of~\eqref{th3.3} and~\eqref{th4.02}  $(-1)^s\Delta(0) < 0$, and hence in virtue of the Lemma~\ref{l::signChange} we get $\Delta^{(s)}(0) < 0$. However this contradicts to~\eqref{th4.2}.

In virtue of proved above property~\eqref{ekstr} the limit $\lim\limits_{M_k\to\infty}\left|G^{(k)}_{r,n,f,M_k}(0) \right|$ does not depend on the choice of the sequence $\{M_k\}_{k=1}^\infty$. This finishes the proof of the theorem.

\subsection{Proof of Theorem~\ref{th::6}.} 

We will need the following lemmas. 

\begin{lemma}\label{l::4}
Let $\varphi_{r,n,f}(\infty) < \infty$ and $Q_{r,n}\in L^r_{\infty,\infty}$ is $g$-spline of $r$-th order defined on half-line $[0,\infty)$ with $n\in{\mathbb Z}_+$ knots $0=:t_0<t_1<\ldots<t_n$. Then there exists $\epsilon\in\{-1,1\}$ such that for $s=1,2,\dots,r$
	\begin{equation}\label{l4.1}
		Q_{r,n}^{(s)}(t)=\epsilon\,P_r^{(s)}(t),\,\,t\geq t_n.
	\end{equation}
\end{lemma}
We will prove lemma using induction. In the case $s=r$ equality~\eqref{l4.1} holds.
Let it be true for $s=k\geq 2$. We will prove, that it is true for $s=k-1$ too. In virtue of the induction assumption $Q_{r,n}^{(k)}(t)=\epsilon\,P_r^{(k)}(t),\,\,t\geq t_n.$ Moreover, $Q_{r,n}^{(k-1)}(\infty)=P_{r}^{(k-1)}(\infty)=0$. Then for $t\geq t_n$  $$-Q_{r,n}^{(k-1)}(t)=\int\limits_{t}^\infty Q_{r,n}^{(k)}(s)ds=\epsilon\int\limits_{t}^\infty P_{r}^{(k)}(s)ds= -\epsilon P_{r}^{(k-1)}(t).$$ Lemma is proved.

\begin{lemma}\label{l::5}
Let $\varphi_{r,n,f}(\infty) < \infty$ and $\lim\limits_{t\to\infty}\frac{f(t)}{P_r(t)} = \infty$. Then the number of oscillation points of $g$-spline $P_{r,n,f}$ tends to infinity when $n\to\infty$. 
\end{lemma}

Let some $n\in {\mathbb N}$ be fixed. Let $M>0$ is such, that for all $t>M$ $\frac{f(t)}{P_r(t)} > \frac{2}{\varphi_{r,n,f}(\infty)}$. Let $g$-spline $P_{r,n,f}$ have $k$ oscillating points $0\leq a_1<a_2<\ldots<a_{k}$.  Denote by $0\leq b_1<b_2<\ldots<b_{n+r+1}$ all oscillation points of $g$-spline $G_{r,n,f,K}$, where $K$ is chosen so big, that ${\rm sgn} G_{r,n,f,K}(b_s)={\rm sgn} P_{r,n,f}(b_s)$, $s=1,2,\dots,k$, $b_{k+1}>\max\{a_{k}, M\}$ and $\varphi_{r,n,f}(K) > \frac 12 \varphi_{r,n,f}(\infty)$.

Set $\Delta(t):= P_{r,n,f}(t)- G_{r,n,f,K}(t)$. Then ${\rm sgn} \Delta(a_s)={\rm sgn} P_{r,n,f}(a_s)$, $s=1,2,\dots,k$. Hence the function $\Delta(t)$ has $k-1$ zeros on the interval $[0,a_{k}]$. Moreover, for $s = k+1,\dots, n+r+1$ $$|P_{r,n,f}(b_s)|< \frac{\varphi_{r,n,f}(\infty)}{2}f(b_s) < \varphi_{r,n,f}(K)f(b_s) = |G_{r,n,f,K}(b_s)|.$$ This means, that ${\rm sgn} \Delta(b_s)={\rm sgn} G_{r,n,f,K}(b_s)$, $s=k+1,\dots,n+r+1$. Hence function $\Delta(t)$ has $n+r-k$ zeros on the interval $[b_{k+1},a]$, and hence on the whole interval $[0,K]$ --- at least $n+r-1$ zeros. This means that function $\Delta^{(r)}(t)$  has at least $n-1$ sign changes and hence $g$-spline $P_{r,n,f}$ has at least $n-1$ knots. 

Let $t_{n,s}^K$, $s = 1,\dots, n$ be the knots of $g$-spline $G_{r,n,f,K}$. Note, that for all $s = 1,2,\dots, n$ on the interval $[0,t_{n,s}^k]$  $g$-spline $G_{r,n,f,K}$ has at least $s$ oscillating points. Really, assume contrary, let for some $1\leq s\leq n$ on interval $[0,t_{n,s}^K]$  $g$-spline $G_{r,n,f,K}$ has less than $s$ oscillation points. Then on the interval $(t_{n,s}^K, K]$ $g$-spline $G_{r,n,f,K}$ has more than $n+r+1-s$ oscillation points, and hence more than $n+r-s$ sign changes. This means that the function $G^{(r)}_{r,n,f,K}$ has more than $n-s$ sign changes on the interval $(t_{n,s}^K, K]$, which is impossible.

This means, that limiting $g$-spline $P_{r,n,f}$ has at least $n-1$ oscillation points. Lemma is proved.

\begin{lemma}\label{l::6}
For all $s=1,2,\dots,r$ and all $t\geq 0$ (almost everywhere is the case when $s=r$) following inequality holds
	\begin{equation}\label{l6.1}
		\left|P^{(s)}_{r,n,f}(t)\right|\leq\left|P^{(s)}_r(t)\right|.
	\end{equation}
\end{lemma}
We will prove the statement of the lemma using induction. In the case $s=r$ inequality~\eqref{l6.1} holds with equality sign. Let inequality~\eqref{l6.1} hold with $s=k\geq 2$. We will prove that it is true for $s=k-1$.

Assume contrary. Let $$T:=\left\{t\in[0,\infty):\left|P^{(k-1)}_{r,n,f}(t)\right|>\left|P^{(k-1)}_r(t)\right|\right\}\neq\emptyset.$$
Denote by $0<t_1<\ldots<t_l$ all knots of $g$-spline $P_{r,n,f}$. Then in virtue of the Lemma~\ref{l::4} $T\subset[0,t_l)$ and 
\begin{equation}\label{l6.2}
\left|P^{(k-1)}_{r,n,f}(t_l)\right|=\left|P^{(k-1)}_r(t_l)\right|.
\end{equation}

Let $a\in T$. Then $$\left|P^{(k-1)}_r(t_l)-P^{(k-1)}_r(a)\right|=\left|\int\limits_a^{t_l}P_r^{(k)}(t)dt\right|=
\int\limits_a^{t_l}\left|P_r^{(k)}(t)\right|dt\geq\int\limits_a^{t_l}\left|P^{(k)}_{r,n,f}(t)\right|dt\geq$$ $$\geq \left|\int\limits_a^{t_l}P^{(k)}_{r,n,f}(t)dt\right|=\left|P^{(k-1)}_{r,n,f}(t_l)-P^{(k-1)}_{r,n,f}(a)\right|,$$
what is impossible in virtue of~\eqref{l6.2}, and the facts that ${\rm sgn} P_r^{(k-1)}(t_l)={\rm sgn} P_r^{(k-1)}(a)$, functions $\left|P_r^{(k-1)}\right|$ is non-increasing and $a\in T$. Contradiction. Hence $T=\emptyset$ and lemma is proved.

Let's return to the proof of the theorem.

Let 
\begin{equation}\label{th01.1}
	\varliminf\limits_{t\to\infty}\frac{f(t)}{\left|P_r(t)\right|}=:c<\infty.
\end{equation} 
Then in virtue of the Theorem~\ref{th::2}  $c>0$. We will show that $\varphi_{r,n,f}(\infty)\geq\frac 1{2c}\:\forall n$. Assume contrary, let a number $n_0$ such, that
\begin{equation}\label{th01.2}
	\varphi_{r,n_0,f}(\infty)<\frac 1{2c}
\end{equation} 
exist. Denote by $0<t_1<t_2<\ldots<t_k$ all the knots of $P_{r,n_0,f}$. Then in virtue of Lemma~\ref{l::4} $P^{\prime}_{r,n_0,f}(t)=\pm\, P_r^{\prime}(t)$, $t\geq t_k$.
In virtue of~\eqref{th01.2} we have $$\left|P_{r,n_0,f}(t)\right|<\frac{f(t)}{2c}$$ for all $t\geq 0$. $f(\infty)=0$, since~\eqref{th01.1} holds. Then $P_{r,n_0,f}(\infty)=0$, and hence $$P_{r,n_0,f}(t)=\pm\, P_r(t),$$ $t\geq t_k$. But then $|P_r(t)|<\frac{f(t)}{2c},$ i. e. $$\frac {f(t)}{|P_r(t)|}>2c,\,t\geq t_n,$$ which contradicts to~\eqref{th01.1}. Sufficiency is proved.

We will prove the necessity now. $\lim\limits_{n\to\infty} \varphi_{r,n,f}(\infty)=\delta>0$. Assume contrary, let $\varliminf\limits_{t\to\infty}\frac{f(t)}{\left|P_r(t)\right|\vphantom{^2}}=\infty$. Then there exists a number $M>0$ such, that for all $t>M$ inequality $$\frac{f(t)}{|P_r(t)|}>\frac 1{\delta}$$ holds, which is equivalent to $$|P_r(t)|<\delta f(t).$$ In virtue of Lemma~\ref{l::6} for all $n$ the following inequality holds
\begin{equation}\label{th01.3}
\int\limits_0^M\left|P^{\prime}_{r,n,f}(t)\right|dt\leq \int\limits_0^M\left|P_{r-1}(t)\right|dt.
\end{equation}

Choose $n$ so big, that 
\begin{equation}\label{th01.4}
	n\cdot f(M)>\int\limits_0^M\left|P_{r-1}(t)\right|dt.
\end{equation}
Choose $m$ such that $g$-spline $P_{r,m,f}(t)$ has at least $n + 1 $ oscillation points (this is possible in virtue of Lemma~\ref{l::5}). Denote by $0\leq a_1<a_2<\ldots<a_{n+1}$ the first oscillation points of $g$-spline $P_{r,m,f}(t)$. Then in virtue of~\eqref{th01.4}, and the facts that $f$ is non-increasing function, $\bigvee\limits_0^M P_{r,m,f}=\int\limits_0^M\left|P^{\prime}_{r,m,f}(t)\right|dt $ and \eqref{th01.3} we get that $a_{n}>M$. Thus $a_{n+1}>a_{n}>M$ are oscillation points of $g$-spline $P_{r,m,f}(t)$ and we get $$\left|P_{r,m,f}(a_{n+1})-P_{r,m,f}(a_{n})\right|=\varphi_{r,m,f}(\infty)\cdot(f(a_{n+1})+f(a_{n}))\geq$$ $$\geq\delta\cdot (f(a_{n+1})+f(a_{n}))>\delta f(a_{n})>|P_r(a_{n})|=\int\limits_{a_{n}}^{\infty}\left|P_{r-1}(t)\right|dt.$$

From the other side in virtue of Lemma~\ref{l::6}
$$\left|P_{r,m,f}(a_{n+1})-P_{r,m,f}(a_{n})\right|\leq\int\limits_{a_{n}}^{a_{n+1}}\left|P^{\prime}_{r,m,f}
(t)\right|dt\leq$$
$$\leq\int\limits_{a_{n}}^{a_{n+1}}\left|P_{r-1}(t)\right|dt<\int\limits_{a_{n}}^{\infty}\left|P_{r-1}(t)\right|dt.$$
Contradiction. Theorem is proved.

\bibliographystyle{plain}
\bibliography{KolmogorovInequalities}  

\end{document}